\documentclass[conference]{IEEEtran}

\IEEEoverridecommandlockouts                              % This command is only needed if % you want to use the \thanks command
\makeatletter
	\let\NAT@parse\undefined
\makeatother

	\PassOptionsToPackage{capitalize,nameinlink}{cleveref}
	\PassOptionsToPackage{noend}{algpseudocode}
	\usepackage[cal=boondoxo]{mathalfa}
	\usepackage{acronym}
	\usepackage[space]{cite}
	\usepackage{textcomp}
	
	\usepackage[
		alglinenumber,
		eqreset=section,
		thmreset=section,
		noload={xpatch,todonotes},
		arxiv,
	]{myPreamble}
	\pgfplotsset{compat=1.16}

	\crefname{enumeratpropi}{property}{properties}
	\crefname{enumeratpropii}{property}{properties}
	\crefname{section}{Section}{Sections}
	\crefname{ALG@line}{Step}{Steps}
	\crefname{figure}{Figure}{Figures}
	\Crefname{figure}{Figure}{Figures}

	\makeatletter
		\renewcommand{\operator@font}{\rm}
	\makeatother
	\DeclareMathOperator{\cont}{C}

	\renewcommand{\l}{\ell}
	\renewcommand{\j}{\mathcal j}

	\newcommand{\danger}{\smash{\raisebox{0.05cm}{{\small\fontencoding{U}\fontfamily{futs}\selectfont\char 66\relax}}}}
	\newcommand{\finger}{{\fontencoding{U}\fontfamily{futs}\selectfont t\relax}}
	\newcommand{\myvec}[1]{\mathbf{#1}} % denote vectors as...
	\newcommand{\myop}[1]{#1} % denote operators as...
	\newcommand{\ite}{{k}} % iteration i
	\newcommand{\iite}{{k+1}} % iteration i+1
	\newcommand{\pite}{{k-1}} % iteration i-1	
	\newcommand{\varx}{\myvec{x}} % variable 1
	\newcommand{\dataf}{g} % data-fitting term
	\newcommand{\reg}{h} % regularizer term
	\newcommand{\dimn}{n} % dimension 2
	\newcommand{\mreauprm}{\lambda}
	\newcommand{\llprm}{\mu}
	\newcommand{\norm}[2]{\left\|{}#1{}\right\|_{#2}} % norm

	\def\xmax{1.75}%
	\def\ymax{25}%

	\pgfplotsset{
		h/.style = {%
			line width = 1.0pt,
		},
		myaxis/.style = {%
			no markers,
			xmin    = -\xmax,
			xmax    = \xmax,
			ymin    = -0.25,
			ymax    = \ymax+0.25,
			restrict y to domain=0:\ymax,
			xlabel  = {},
			ylabel  = {},
			title   = {},
			xtick        = \empty,
			ytick        = \empty,
			axis lines   = middle,
			axis line style = {draw=none},
			legend style = {%
				rounded corners,
				fill opacity      = 0.7,
				font              = \color{black}\scriptsize,
				legend cell align = left,
				anchor            = south west,
				at                = {(0,0)},
				line width        = 0.4pt,
				inner sep         = 1pt,
				outer sep         = 5pt,
				cells             = {line width=1.5pt},
			},
			execute at begin axis = {
				\node[inner sep=0pt, outer sep=0pt] at (axis cs:0,-0.25) {};
				\node[inner sep=0pt, outer sep=0pt] at (axis cs:0,\ymax+0.25) {};
			},
		},
	}

	\tikzset{%
		declare function = {%
			h(\x)         = (abs(abs(\x)-1) > 0) * \ymax;                        % h = indicator of {1,-1}
			hM(\x,\l)     = (abs(\x)-1)^2 / (2*\l);                              % Moreau envelope of h
			hHull(\x,\l)  = (abs(\x)<= 1) * (1-\x^2)/(2*\l) + (abs(\x)>1)*\ymax; % proximal hull of h
			Ph(\x,\l,\m)  = (abs(\x) > \m/\l) * (\l*\x-sign(\x)*\m) / (\l-\m)  ; % argmax in the LL subproblem
			hLL(\x,\l,\m) = hM(Ph(\x,\l,\m), \l) - (Ph(\x,\l,\m)-\x)^2/(2*\m)  ; % LL envelope of h
		},
		Dot/.style = {%
			circle,
			fill,
			draw = black,
			inner sep = 1.5pt,
			outer sep = 0pt,
		},
		Void/.style = {%
			rectangle,
			fill = white,
			draw = none,
			inner sep = 3pt,
			outer sep = 0pt,
		},
	}

\showgrayfalse
\begin{document}

\title{%
	Lasry--Lions Envelopes and Nonconvex Optimization: A Homotopy Approach\thanks{%
		This work was supported by the Research Foundation -- Flanders (FWO) projects G0A0920N, G086518N, and G086318N, by the Research Council KU Leuven C1 project C14/18/068, and by the Fund for Scientific Research -- FNRS and FWO EOS project 30468160 (SeLMA).%
	}%
}

\author{\IEEEauthorblockN{Miguel Sim\~{o}es}
\IEEEauthorblockA{\textit{Dept.\ Electr.\ Eng.\ (ESAT) -- STADIUS} \\
\textit{KU Leuven}\\
Leuven, Belgium \\
miguel.simoes@kuleuven.be}
\and
\IEEEauthorblockN{Andreas Themelis}
\IEEEauthorblockA{\textit{Fac.\ Inf.\ Sci.\ Electr.\ Eng.\ (ISEE)} \\
\textit{Kyushu University}\\
Fukuoka, Japan \\
andreas.themelis@ees.kyushu-u.ac.jp}
\and
\IEEEauthorblockN{Panagiotis Patrinos}
\IEEEauthorblockA{\textit{Dept.\ Electr.\ Eng.\ (ESAT) -- STADIUS} \\
\textit{KU Leuven}\\
Leuven, Belgium \\
panos.patrinos@kuleuven.be}%
}

\maketitle

	\begin{abstract}
		In large-scale optimization, the presence of nonsmooth and nonconvex terms in a given problem typically makes it hard to solve.
		A popular approach to address nonsmooth terms in convex optimization is to approximate them with their respective Moreau envelopes.
		In this work, we study the use of Lasry--Lions double envelopes to approximate nonsmooth terms that are also not convex.
		These envelopes are an extension of the Moreau ones but exhibit an additional smoothness property that makes them amenable to fast optimization algorithms. Lasry--Lions envelopes can also be seen as an ``intermediate'' between a given function and its convex envelope, and we make use of this property to develop a method that builds a sequence of approximate subproblems that are easier to solve than the original problem. We discuss convergence properties of this method when used to address composite minimization problems; additionally, based on a number of experiments, we discuss settings where it may be more useful than classical alternatives in two domains: signal decoding and spectral unmixing.
	\end{abstract}

	%taken from https://www.ieee.org/content/dam/ieee-org/ieee/web/org/pubs/taxonomy_v101.pdf

	\section{Introduction}
		Many problems in signal processing are formulated as composite problems of the form
\begin{equation} \label{eq:composite_prob}
\minimize_{\varx \in \R^\dimn} \varphi(\varx) \coloneqq \dataf(\varx) +  \reg(\varx),
\end{equation}
where $\dataf:\R^\dimn \to \Rinf$ is a data-fitting term and $\reg:\R^\dimn \to \Rinf$ acts as a regularizer or as a constraint.
These problems are typically faced with a number of challenges: they involve a very large number of variables, they are ill posed,
%---in the sense that they do not satisfy any of the three conditions suggested by Hadamard for a problem to be considered well posed (stability, existence, and uniqueness of the solution)---
and their data-fitting and regularizer terms are nonsmooth and nonconvex.
As an example, consider that one wishes to solve a regression problem whose solution is known to be sparse.
One may use the so-called $\ell_0$ pseudo-norm as a regularizer, i.e., $\reg = \norm{\cdot}{0}$, which measures the number of nonzeros entries of a given vector.
Finding global solutions to large-scale problems involving it is impractical, since this function is nonsmooth and nonconvex, and these problems are often approached by replacing it with surrogates, such as the $\ell_1$ norm, which allow one to use methods that are more computationally tractable.
Under certain conditions on $\dataf$, both approaches produce equivalent solutions~\cite{Candes2011}, but these conditions are often violated.
Consequently, it is of interest to solve the problem with the nonconvex regularizer.
Perhaps surprisingly, in practice it is not always the case that the use of the $\ell_0$ pseudo-norm instead of the $\ell_1$ surrogate will produce a solution that is more useful, due to either the inherent ill posedness of the problem or wrong modeling assumptions on $\dataf$.
However, recent work has presented experimental evidence that this solution may in fact be more useful in certain noise regimes~\cite{Bertsimas2016,Hastie2020}, and has validated the idea that a combination of the $\ell_0$ regularizer with additional $\ell_1$ or $\ell_2$ terms produces the most useful results for the different noise regimes of practical interest~\cite{Hazimeh2020}.

In this work, we discuss a method that makes use of the Lasry--Lions double envelope to build surrogates of nonsmooth nonconvex functions.
For a proper, lower semicontinuous (lsc) function $\reg:\R^\dimn\to\Rinf$ and parameters $\mreauprm>\llprm>0$, its $\mreauprm$-\emph{Moreau envelope} is
\begin{align*}
	\reg^{\mreauprm}(\myvec{x}) \coloneqq \inf_\myvec{w}\left\{\reg(\myvec{w})+\tfrac{1}{2\mreauprm} \norm{ \myvec{w}-\myvec{x}}{}^2\right\}
\end{align*}
and its $\left(\mreauprm,\llprm\right)$-\emph{Lasry--Lions double envelope} is 
\begin{align}\label{eq:doublenv}
	\reg^{\mreauprm,\llprm}(\myvec{x}) \coloneqq -(-\reg^{\mreauprm})^\llprm (\myvec{x}) = \sup_\myvec{w}\left\{\reg^{\mreauprm}(\myvec{w})-\tfrac{1}{2\llprm} \norm{\myvec{w}-\myvec{x}}{}^2\right\}.
\end{align}
Replacing $\reg$ in~\eqref{eq:composite_prob} with its envelope $\reg^{\mreauprm,\llprm}$---which, e.g., can be made for the $\ell_0$ pseudo-norm or the indicator of a nonconvex set---furnishes approximate subproblems arbitrarily close to the original one and that are, in a way to be made clearer in a later section, ``less'' nonconvex.
We also show that these problems are smooth, and hence can be tackled via fast smooth optimization solvers.
We make use of the surrogate problems in a way that mimics the behavior of homotopy methods for systems of nonlinear equations or of penalty methods for constrained optimization: we solve a sequence of subproblems that increasingly approximate the original problem, and use the solution of a given subproblem to inform the initial estimate of the following one.
As we show in \cref{sec:min_algos}, the sequence of global minimizers of these subproblems converges to the global minimizer of~\eqref{eq:composite_prob}.
We motivate the use of this method in two ways.
First, we argue that the use of Lasry--Lions double envelopes provides a way to systematically produce surrogates of nonconvex functions that is useful from a computational perspective, in the sense just described.
Second, its use also seems to have benefits in terms of the solutions that are found through it.
In fact, the Lasry--Lions envelope of the $\ell_0$ pseudo-norm can be seen as a combination of nonconvex and convex terms. This combination is different from the combinations of the $\ell_0$ regularizer with convex terms described above but, when compared to the use of the pure $\ell_0$ regularizer, can be said to be less ``aggressive'', in the sense used in~\cite{Hastie2020} when referring to the bias--variance trade-off issue present in estimation problems.

To the best of our knowledge, the use of the Lasry--Lions envelope as described seems to be new, but there are some connections to existing approaches in the literature.
For example, there is a family of methods that addresses the problems of $\ell_0$-regularized least squares and of rank-constrained matrix approximation by replacing the nonconvex terms with others that are still nonconvex but more tractable; additionally, results associating the minimizers of both formulations exist~\cite{Soubies2017,Carlsson2019,Grussler2018}.
The connections to our work come from the fact that these replacements correspond to the proximal hull of certain functions---the {\(\mreauprm\)-proximal hull} of function $\reg$ corresponds to making \(\llprm=\mreauprm\) in~\eqref{eq:doublenv}, which we denote by \(\reg^{\mreauprm,\mreauprm}\).
Unlike the Lasry--Lions envelope, the proximal hull is not necessarily smooth.

The structure of this work is as follows.
In \cref{sec:double_env}, we briefly discuss a number of properties of the Lasry--Lions double envelope as well as some examples of interest.
In \cref{sec:min_algos}, we study the homotopy approach that makes use of the Lasry--Lions envelope, and discuss a practical implementation of it.
In \cref{sec:numer_exps}, we present numerical experiments in signal decoding and spectral unmixing.
\Cref{sec:conclusions} concludes.

	\section{The Lasry--Lions double envelope} \label{sec:double_env}
		We start by stating our assumptions for~\eqref{eq:composite_prob} and then list some results concerning the Lasry--Lions double envelope: \cref{thm:LL=,thm:LL_prop} are concerned with some of its basic identities and properties, and \cref{thm:LLsmooth,thm:LLgrad:equiv} with its smoothness characteristics.
We finish the section with two examples of double envelopes: the indicator function of the set $\{-1,1\}$ in \cref{thm:ex_ind__binary} and the $\ell_0$ pseudo-norm in \cref{thm:ex_l0}; we illustrate the latter in \cref{fig:LL_ind_binary}.

\begin{figure*}[htbp]
		\centering
	
\begin{minipage}[t]{0.3\linewidth}%
	\includetikz[width=\linewidth]{biconj}
	\hspace*{0pt}\hfill{\it From the convex relaxation}\hfill\hspace*{0pt}
	
	\vspace{3pt}\small
	As \(\lambda>\mu\to\infty\) with \(\lambda-\mu\to0\),
% 	the Lasry--Lions envelope
	\(h^{\lambda,\mu}\) converges to the convex hull of the original function \(h^{**}=\indicator_{[-1,1]}\) (red).%
\end{minipage}%
\hfill
\begin{minipage}[t]{0.3\linewidth}%
	\includetikz[width=\linewidth]{orig}
	\hspace*{0pt}\hfill{\it to the original problem}\hfill\hspace*{0pt}

	\vspace{3pt}\small
	As \(0<\mu<\lambda\to0\),
% 	the Lasry--Lions envelope
	\(h^{\lambda,\mu}\) converges to the original function \(h=\indicator_{\{\pm1\}}\).%
\end{minipage}%
\hfill
\begin{minipage}[t]{0.3\linewidth}%
	\includetikz[width=\linewidth]{Moreau}
	\hspace*{0pt}\hfill{\it in a smooth fashion.}\hfill\hspace*{0pt}
	
	\vspace{3pt}\small
	Approximations of extended-real proximal hulls,
	continuous Moreau envelopes,
	and smooth Lasry--Lions envelopes.%
\end{minipage}%	
	
	\caption{Approximations with Lasry--Lions envelopes \(h^{\lambda,\mu}\) (blue), proximal hulls \(h^{\lambda,\lambda}\) (brown), and Moreau envelopes \(h^\lambda\) (green) of the nonconvex, nonsmooth, extended-real-valued function \(h=\indicator_{\{\pm1\}}\) (black).}
	\label{fig:LL_ind_binary}
\end{figure*}

A detailed account of the notions used throughout this work can be found in~\cite{Rockafellar2009}.
%We denote the \emph{scalar product} of a Hilbert space by $\langle \cdot , \cdot \rangle$ and the associated \emph{norm} by $\| \cdot \|$.
% The Legendre--Fenchel \emph{conjugate} of a function $f$ is denoted by $f^*$.
% The \emph{indicator function} of a set $C \in \R^{n}$ is defined as $\delta_C(\mathbf{x}) \defeq 0$ if $\mathbf{x} \in C$, $\delta_C(\mathbf{x}) \defeq + \infty$ otherwise.
Bold lowercase letters denote vectors, bold uppercase letters denote matrices, $[\mathbf{a}]_i$ denotes the $i$-th element of a vector $\mathbf{a}$, and $[\mathbf{A}]_{:j}$ denotes the $j$-th column of a matrix $\mathbf{A}$.
%, and $[\mathbf{A}]_{ij}$ denotes the element in the $i$-th row and $j$-th column of a matrix $\mathbf{A}$. 
%The signum operator is denoted by $\text{sgn}(\cdot)$.
We let \(\j\coloneqq\tfrac12\|{}\cdot{}\|^2\).

\begin{ass}\label{ass:basic}%
	In problem~\eqref{eq:composite_prob}, the following hold:
	\begin{enumeratass}
		\item\label{ass:g}%
		\(\dataf\in\cont^1(\R^n)\);
		\item\label{ass:h}%
		\(\func{\reg}{\R^m}{\Rinf}\) is proper, lsc, and prox-bounded;
		%\label{ass:F}%
		%\(F\in\cont^1(\R^n;\R^m)\) with \(\J F\) surjective
		\item\label{ass:phi}%
		a solution exists: \(\argmin\varphi\neq\emptyset\).
	\end{enumeratass}
\end{ass}

\begin{fact}[Basic identities of the Lasry--Lions envelope]\label{thm:LL=}%
	The following hold for \(\func{h}{\R^n}{\Rinf}\) and \(\lambda\geq\mu>0\):
	\begin{enumerate}
		\item\label{thm:moreau_proxh_double}%
		\(
		h^{\lambda,\mu}
		{}={}
		(h^{\lambda,\lambda})^{\lambda-\mu}
		{}={}
		(h^{\lambda-\mu})^{\mu,\mu}
		\);
		\item\label{thm:MoreauConj}%
		\(
		h^\lambda
		{}={}
		\lambda^{-1}\j
		{}-{}
		\conj{\bigl(
			h+\lambda^{-1}\j
			\bigr)}(\lambda^{-1}{}\cdot{})
		\);
		\item\label{thm:proxh}%
		\(
		h^{\lambda,\lambda}
		{}={}
		\biconj{\bigl(
			h+\lambda^{-1}\j
			\bigr)}
		{}-{}
		\lambda^{-1}\j
		\).
	\end{enumerate}
	When \(\lambda>\mu\), denoting \(c\coloneqq\frac{\lambda(\lambda-\mu)}{\mu}\) it also holds that%
	\begin{enumerate}[resume]
		\item\label{thm:LL=conj}%
		\(
		h^{\lambda,\mu}
		{}={}
		(\lambda-\mu)^{-1}\j
		{}-{}
		\bigl[
		\conj{\bigl(
			h+\lambda^{-1}\j
			\bigr)}
		\bigr]^{\nicefrac1c}
		\bigl(
		\tfrac{{}\cdot{}}{\lambda-\mu}
		\bigr)
		\);
		\item\label{thm:LL=conj2}%
		\(
		h^{\lambda,\mu}
		{}={}
		\bigl[
		\biconj{\bigl(
			h+\lambda^{-1}\j
			\bigr)}
		\bigr]^c
		\bigl(
		\tfrac\lambda\mu{}\cdot{}
		\bigr)
		{}-{}
		\mu^{-1}\j
		\).
	\end{enumerate}
	\begin{proof}[Detail]
		\ref{thm:moreau_proxh_double}--\ref{thm:proxh}: follow from~\cite[Ex.s 1.46, 11.26(c)]{Rockafellar2009}; \ref{thm:LL=conj}, \ref{thm:LL=conj2}: from \ref{thm:moreau_proxh_double}--\ref{thm:proxh}.
	\end{proof}
\end{fact}

\begin{fact}[Basic properties of the Lasry--Lions envelope] \label{thm:LL_prop}%
	Let \(\func{h}{\R^n}{\Rinf}\) be proper, lsc, and \(\gamma_h\)-prox-bounded.
	Then, \(h^{\lambda,\mu}\) is proper and lsc for every \(0<\mu\leq\lambda<\gamma_h\).
	Moreover%
	\begin{enumerate}
	\item\label{thm:sandw}
		\(
		h^\lambda
		{}\leq{}
		h^{\lambda,\mu}
		{}\leq{}
		h^{\lambda-\mu}
		{}\leq{}
		h
		\);
	\item\label{thm:approx}
		\(h^{\lambda,\mu}\) is pointwise increasing wrt \(\mu\) and decreasing both wrt \(\lambda\) and \(\lambda-\mu\)
		(\ie, \(h^{\lambda,\mu}\geq h^{\lambda',\mu'}\) when \(\lambda-\mu\leq\lambda'-\mu'\));
	\item\label{thm:LL:pointwise}%
		for every \(\myvec{\bar x} \in \R^\dimn\), \(h^{\lambda,\mu}(\myvec{x})\to h(\myvec{\bar x})\) as \(\myvec{x}\to \myvec{\bar x}\) and \(0<\mu\leq\lambda\to0\) in such a way that \(\frac{\|\myvec{x}-\myvec{\bar x}\|}{\lambda}\) is bounded;
	\item\label{thm:convdom}
		\(\dom h^{\lambda,\lambda}=\conv\dom h\);
	\item\label{thm:infLL}
		\(\inf h^{\lambda,\mu}=\inf h\) and \(\argmin h^{\lambda,\mu}=\argmin h\);
% 	\item\label{thm:hypocvx}
% 		\(h^{\lambda,\mu}\) is $\mu$-proximal;
	\item\label{thm:LLseparable}
		let \(\func{\reg}{\R^{n_1} \times \cdots \times \R^{n_m}}{\Rinf}\) be of the form \(\reg = \sum_{\ite=1}^{m} \reg_\ite(\myvec{x}_\ite)\), where $\myvec{x}_\ite \in \R^{n_\ite}$ for all $\ite \in \{1,\cdots,m\}$; then
		\[\textstyle
			\reg^{\lambda,\mu}(\myvec{x}) = \sum_{\ite=1}^{m} \reg_\ite^{\lambda,\mu}(\myvec{x}_\ite).
		\]
	\end{enumerate}
	\begin{proof}[Detail]
		\ref{thm:sandw}, \ref{thm:infLL}: follow from~\cite[1.46]{Rockafellar2009}; \ref{thm:approx}--\ref{thm:infLL}: from~\cite[1.25, 1.44]{Rockafellar2009} and~\eqref{thm:LL=}; \ref{thm:LLseparable}: from~\eqref{eq:doublenv}.
	\end{proof}
\end{fact}

\begin{prop}[Smoothness of the Lasry--Lions envelopes]\label{thm:LLsmooth}%
	Let \(\func{h}{\R^n}{\Rinf}\) be proper, lsc, and \(\gamma_h\)-prox-bounded.
	Then, for every \(0<\mu<\lambda<\gamma_h\) the Lasry--Lions envelope \(h^{\lambda,\mu}\) is \(L_{h^{\lambda,\mu}}\)-Lipschitz-continuously differentiable with gradient
	\begin{equation}\label{eq:gradLL}
	\nabla h^{\lambda,\mu}
	{}={}
	\tfrac\lambda\mu
	P
	{}-{}
	\tfrac1\mu\id,
	\end{equation}
	where
	\(
	P
	{}\coloneqq{}
	\prox_{c^{-1}\conj{(h+\lambda^{-1}\j)}}(
	\frac{{}\cdot{}}{\lambda-\mu}
	)
	\)
	and \(c\coloneqq\frac{\lambda(\lambda-\mu)}{\mu}\) is as in \cref{thm:LL=}.
	In fact, for every \(\myvec{x},\myvec{y}\in\R^n\) it holds that
	\[
	\sigma_{h^{\lambda,\mu}}\|\myvec{x}-\myvec{y}\|^2
	{}\leq{}
	\innprod*{
		\nabla h^{\lambda,\mu}(\myvec{x})
		{}-{}
		\nabla h^{\lambda,\mu}(\myvec{y})
	}{\myvec{x}-\myvec{y}}
	{}\leq{}
	-\sigma_{-h^{\lambda,\mu}}\|\myvec{x}-\myvec{y}\|^2
	\]
	for some
	\(
		\sigma_{h^{\lambda,\mu}}\geq\frac{-1}{\mu}
	\)
	and
	\(
	\sigma_{-h^{\lambda,\mu}}\geq\frac{-1}{\lambda-\mu}
	\),
	and in particular
	\(
	L_{h^{\lambda,\mu}}
	{}\leq{}
	\max\set{\frac1\mu,\frac{1}{\lambda-\mu}}
	\).
	Moreover, the estimates can be tightened as
	\begin{enumerate}
		\item\label{thm:LLsmooth:hypo}%
		\(
		\sigma_{h^{\lambda,\mu}}
		{}\geq{}
		\frac{\sigma_h}{1+(\lambda-\mu)\sigma_h}
		\)
		and
		\(
		L_{h^{\lambda,\mu}}
		{}\leq{}
		\max\set{\frac{|\sigma_h|}{1+(\lambda-\mu)\sigma_h},\frac{1}{\lambda-\mu}}
		\)
		if \(h\) is \(\sigma_h\)-hypoconvex;
		\item
		\(
		\sigma_{-h^{\lambda,\mu}}
		{}\geq{}
		-\frac{L_h}{1+(\lambda-\mu)L_h}
		\)
		and
		\(
		L_{h^{\lambda,\mu}}
		{}\leq{}
		\frac{L_h}{1+(\lambda-\mu)L_h}
		\)
		if \(h\) is \(L_h\)-smooth.
	\end{enumerate}
	\begin{proof}
		The expression for the gradient comes from the identity in \cref{thm:LL=conj} together with the known formula \(\nabla\phi^{1/c}=c(\id-\prox_{\phi/c})\), holding for the proper, lsc, convex function \(\phi=\conj{(h+\lambda^{-1}\j)}\).
		In turn, the general bounds on the inner product follow from the monotonicity and nonexpansiveness of \(\prox_{\phi/c}\).
		When \(h\) is \(\sigma_h\)-hypoconvex, then \(\phi\) is \(L_\phi\)-smooth (and \((-L_\phi)\)-hypoconvex) with \(L_\phi=\frac{\lambda}{1+\lambda\sigma_h}\), and consequently \(\prox_{\phi/c}\) is \((1+c^{-1}L_\phi)^{-1}\)-strongly monotone \cite[Lem. 5]{themelis2020new}.
		Similarly, when \(h\) is \(L_h\)-smooth (hence \((-L_f)\)-hypoconvex), then \(\phi\) is \(\sigma_\phi\)-strongly convex with \(\sigma_\phi=\frac{\lambda}{1+\lambda L_h}\), and \(\prox_{\phi/c}\) is thus \((1+c^{-1}\sigma_\phi)\)-contractive.
		By using these estimates, the tighter values as in the statement are obtained.
	\end{proof}
\end{prop}

%\begin{cor}[Smoothness of the Lasry--Lions envelopes when \(\lambda=2\mu\)] \label{thm:smoothness} %
%	Let \(\func{h}{\R^n}{\Rinf}\) be proper, lsc, and \(\gamma_h\)-prox-bounded.
%	Then, for every \(0<\lambda<\gamma_h\) the Lasry--Lions envelope \(h^{\lambda,\lambda/2}\) is \(L_{h^{\lambda,\lambda/2}}\)-Lipschitz-continuously differentiable with gradient
%	\[
%	\nabla h^{\lambda,\lambda/2}
%	{}={}
%	2P
%	{}-{}
%	\tfrac2\lambda\id,
%	\]
%	where
%	\(
%	P
%	{}\coloneqq{}
%	\prox_{\lambda^{-1}\conj{(h+\lambda^{-1}\j)}}(
%	\frac2\lambda {}\cdot{}
%	)
%	\).
%	In fact, for every \(x,y\in\R^n\) it holds that
%	\[
%	\sigma_{h^{\lambda,\lambda/2}}\|x-y\|^2
%	{}\leq{}
%	\innprod*{
%		\nabla h^{\lambda,\lambda/2}(x)
%		{}-{}
%		\nabla h^{\lambda,\lambda/2}(y)
%	}{x-y}
%	{}\leq{}
%	-\sigma_{-h^{\lambda,\lambda/2}}\|x-y\|^2
%	\]
%	for some
%	\(
%	\sigma_{h^{\lambda,\lambda/2}},
%	\sigma_{-h^{\lambda,\lambda/2}}
%	{}\geq
%	-\frac2\lambda
%	\),
%	and in particular
%	\(
%	L_{h^{\lambda,\lambda/2}}
%	{}\leq{}
%	\frac2\lambda
%	\).
%	Moreover, the estimates can be tightened as
%	\begin{enumerate}
%		\item
%		\(
%		\sigma_{h^{\lambda,\lambda/2}}
%		{}\geq{}
%		\frac{2\sigma_h}{2+\lambda\sigma_h}
%		\)
%		if \(h\) is \(\sigma_h\)-hypoconvex,
%		\item
%		\(
%		\sigma_{-h^{\lambda,\lambda/2}}
%		{}\geq{}
%		-\frac{2L_h}{2+\lambda L_h}
%		\)
%		and
%		\(
%		L_{h^{\lambda,\lambda/2}}
%		{}\leq{}
%		\frac{2L_h}{2+\lambda L_h}
%		\)
%		if \(h\) is \(L_h\)-smooth.
%	\end{enumerate}
%\end{cor}

\begin{cor}\label{thm:LLgrad:equiv}%
	For \(\myvec{x}\in\R^n\) and \(0<\mu<\lambda<\gamma_h\) it holds that
	\[
	\myvec{d}=\nabla h^{\lambda,\mu}(\myvec{x})
	\quad\Leftrightarrow\quad
	\myvec{x}-(\lambda-\mu)\myvec{d}
	{}\in{}
	\conv\prox_{\lambda h}(\myvec{x}+\mu \myvec{d}).
	\]
	\begin{proof}
		Denoting \(c\coloneqq\frac{\lambda(\lambda-\mu)}{\mu}\), it follows from \eqref{eq:gradLL} that
		\(
		\myvec{d}=\nabla h^{\lambda,\mu}(x)
		\)
		holds iff
		\begin{align*}
		&\tfrac{\mu \myvec{d}+\myvec{x}}{\lambda}
		{}={}
		\prox_{c^{-1}\conj{(h+\lambda^{-1}\j)}}\bigl(
		\tfrac{\myvec{x}}{\lambda-\mu}
		\bigr) \\
		&\qquad \Leftrightarrow
		c\left(
		\tfrac{\myvec{x}}{\lambda-\mu}
		{}-{}
		\tfrac{\mu \myvec{d}+\myvec{x}}{\lambda}
		\right)
		{}\in{}
		\partial\conj{(h+\lambda^{-1}\j)}\left(
		\tfrac{\mu \myvec{d}+\myvec{x}}{\lambda}
		\right),
		\end{align*}
		where the last equivalence owes to convexity of \(\conj{(h+\lambda^{-1}\j)}\).
		It follows from \cref{thm:MoreauConj} that the latter subdifferential expands to
		\(
		\mu \myvec{d}+\myvec{x}
		{}+{}
		\lambda
		\partial(-h^\lambda)(\mu \myvec{d}+\myvec{x})
		\),
		which in turn equals
		\(
		\conv\prox_{\lambda h}(\myvec{x}+\mu \myvec{d})
		\)
		by virtue of~\cite[Ex. 10.32]{Rockafellar2009}.
		A direct computation reveals that the left-hand side of the inclusion is \(\myvec{x}-(\lambda-\mu)\myvec{d}\).
	\end{proof}
\end{cor}

\begin{es}\label{thm:ex_ind__binary}
	Consider \(h:\R\to\R:x\mapsto \indicator_{\set{\pm1}}(x)\).
	Then, its Moreau envelope and proximal hull are
	\begin{align*}
	h^{\lambda}(x)
	{}={}  %\frac{1}{2\lambda} \min\{(1-x)^2,(1+x)^2)\} {}={}
	\frac{(1-|x|)^2}{2\lambda},
	%\shortintertext{}
	\qquad
	h^{\lambda,\lambda}(x)
	{}={}
	\begin{ifcases}
	+ \infty \;  |x|>\nicefrac\mu\lambda
	\\
	\frac{1-x^2}{2\lambda} \; \text{otherwise},
	\end{ifcases}
	\shortintertext{and its Lasry--Lions envelope and respective gradient are}
	h^{\lambda,\mu}(x)
	{}={} 
	\begin{ifcases}
	\frac{(1-|x|)^2}{2(\lambda-\mu)} \; |x|>\nicefrac\mu\lambda
	\\
	\frac{1}{2\lambda}-\frac{x^2}{2\mu} \; \text{otherwise},
	\end{ifcases}
	\,
	%\shortintertext{and}
	\nabla h^{\lambda,\mu}(x)
	{}={}
	\begin{ifcases}
	\frac{x-\sign x}{\lambda-\mu} \; |x|>\nicefrac\mu\lambda
	\\
	-\frac x\mu \; \text{otherwise}.
	\end{ifcases}
	\end{align*}
	Additionally, from~\cref{thm:LLseparable}, for \(h:\R^\dimn\to\R:\myvec{x}\mapsto \indicator_{\set{\pm1}^\dimn} (\myvec{x})\) we have \(h^{\lambda,\mu}(\myvec{x}) = \sum_{\ite=1}^{n} \indicator_{\set{\pm1}}^{\lambda,\mu} ([\myvec{x}]_\ite)\).
\end{es}

\begin{es} \label{thm:ex_l0}
	Consider $h:\R\to\R:x\mapsto |x|_0=\begin{cases}1&,\quad x\neq 0\\ 0&,\quad x=0\end{cases}$.
	Then,	
	\begin{align*}
		\reg^{\lambda}(x)
	{}={} &
		\begin{ifcases}
			\tfrac{1}{2\lambda}x^2& |x|\leq \sqrt{2\lambda}
		\\
			1& |x|\geq \sqrt{2\lambda}
		\end{ifcases}
	\\
		\reg^{\lambda,\lambda}(x)
	{}={} &
		\begin{ifcases}
			1-\tfrac{1}{2\lambda}\left(|x|-\sqrt{2\lambda}\right)^2
			&
			|x|\leq \sqrt{2\lambda}
		\\
			1&|x|\geq \sqrt{2\lambda}
		\end{ifcases}
	\\
		\reg^{\lambda,\mu}(x)
	{}={} &
		\begin{ifcases}
			\tfrac{1}{2(\lambda-\mu)}x^2
			&
			|x|\leq (1-\nicefrac{\mu}{\lambda})\sqrt{2\lambda}
		\\
			1-\tfrac{1}{2\mu}\left(|x|-\sqrt{2\lambda}\right)^2
			&
			|x|\in[(1-\nicefrac{\mu}{\lambda})\sqrt{2\lambda},\sqrt{2\lambda}]
		\\
			1
			&
			|x|\geq\sqrt{2\lambda}
		\end{ifcases}
	\shortintertext{and}
		\nabla h^{\lambda,\mu}(x)
	{}={} &
		\begin{ifcases}
			\tfrac{1}{(\lambda-\mu)}x
			&
			|x|\leq (1-\nicefrac{\mu}{\lambda})\sqrt{2\lambda}
		\\
			-\tfrac{1}{\mu}\left(|x|-\sqrt{2\lambda}\right)
			&
			|x|\in[(1-\nicefrac{\mu}{\lambda})\sqrt{2\lambda},\sqrt{2\lambda}]
		\\
			0
			\otherwise[\(
				|x|\geq \sqrt{2\lambda}.
			\)]
		\end{ifcases}
	\end{align*}
	For \(h:\R^\dimn\to\R:\myvec{x}\mapsto \norm{(\myvec{x})}{0}\), \(h^{\lambda,\mu}(\myvec{x}) = \sum_{\ite=1}^{n} \left(|[\myvec{x}]_\ite|_0\right)^{\lambda,\mu}\).
\end{es}

	\section{A homotopy approach to nonconvex minimization} \label{sec:min_algos}
		This section describes a new method to address nonsmooth nonconvex composite problems.
We state global convergence properties and then discuss a practical implementation of it.

The approximating properties of the Lasry--Lions envelope naturally lead one to consider homotopy approaches to address~\eqref{eq:composite_prob}.
When \cref{ass:h} is satisfied, one has that
\[
\varphi(\myvec{x})
{}={}
\lim_{\mu<\lambda\to0}\set{
	g(\myvec{x})+h^{\lambda,\mu}(\myvec{x})
}
\]
holds for every \(\myvec{x}\in\R^n\), and \(h^{\lambda,\mu}\) is a Lipschitz-differentiable function whenever \(0<\mu<\lambda<\gamma_h\).
In what follows, it will be assumed without mention that \(\lambda\) and \(\mu\) comply with these bounds.
By replacing the regularizer in~\eqref{eq:composite_prob} with its Lasry--Lions envelope, we have:
\[\label{eq:Prs}
\minimize_{\myvec{x}\in\R^n} \varphi_{\lambda,\mu}(\myvec{x})\coloneqq g(\myvec{x})+h^{\lambda,\mu}(\myvec{x}).
\]

\begin{thm}\label{thm:subseq}%
	Additionally to \cref{ass:basic}, suppose that there exist \(0<\bar \mu<\bar \lambda<\gamma_h\) such that \(\varphi_{\bar \lambda, \bar \mu}\) is level bounded.
	Let \(\lambda_k,\mu_k>0\) satisfy
	\(
		\bar\lambda
	{}\geq{}
		\lambda_k
	{}>{}
		\mu_k
	{}\leq{}
		\bar\mu
	\)
	with \(\lambda_k\to0\) and \(\lambda_k-\mu_k\searrow0\) (so that \(\varphi_{\lambda_k,\mu_k}\nearrow\varphi\)).
	The following hold:%
	\begin{enumerate}
		\item\label{thm:penalty:inf}%
			\(\min\varphi_{\lambda_k,\mu_k}\to\min\varphi\) as \(k\to\infty\);
		\item\label{thm:penalty:argmin}%
			any sequence \(\myvec x_\star^k\in\varepsilon_k\)-\(\argmin\varphi_{\lambda_k,\mu_k}\) with \(\varepsilon_k\to0\) is bounded and has all its cluster points in \(\argmin\varphi\).
	\end{enumerate}
	\begin{proof}
		\cref{thm:approx,thm:LL:pointwise} ensure through \cite[Prop. 7.4(d)]{Rockafellar2009} that \(\varphi_{\lambda_k,\mu_k}\) epi-converges to \(\varphi\) as \(k\to\infty\).
		In turn, owing to level boundedness the proof follows by invoking \cite[Ex. 7.32(a) and Thm. 7.33]{Rockafellar2009}.
	\end{proof}
\end{thm}

A conceptual method informed by these results is given in \cref{algo:meta_smooth}.
For ease of implementation, the \(\varepsilon_k\)-minimality prescribed in \cref{thm:subseq} is here replaced by approximate stationarity, so that any smooth minimization algorithm can conveniently be employed.
Extending~\cref{thm:subseq} to account for this relaxation is planned for future work, possibly by mimicking arguments from~\cite{birgin2014practical,grapiglia2020complexity,evens2021neural} made in the context of augmented Lagrangian methods.

In more detail, for the experiments detailed in \cref{sec:numer_exps}, we made the following choices: we let  $\tau^1 = 10^{-3}$; \cref{algo:meta_mininizer} can be implemented by any smooth minimization algorithm, and we used L-BFGS due to its fast properties and low memory requirements; finally, we updated the tolerance as $\tau^\iite = 0.9 \tau^\ite$.

\begin{algorithm}[tbp]
% 	\caption{Lasry--Lions homotopy}%
	\caption{A homotopy approach to address~\eqref{eq:composite_prob} through the use of Lasry--Lions envelopes.}%
	\label{algo:meta_smooth}%
	% \begin{algorithmic}[1]
% 	\Statex Choose $\myvec{x}^{0,+} \in \R^\dimn$, $\lambda^1 > \mu^1 > 0$ and $\tau^1 > 0$ \label{algo:initial_outer}
% 	\Statex $\ite \leftarrow 1$
% 	\While{stopping criterion is not satisfied}
% 	\State Find an approximate minimizer $\myvec{x}^{\ite,+}$ of $\varphi_{\lambda_k,\mu_k}(\cdot)$ with starting point $\myvec{x}^{\ite,1} \coloneqq \myvec{x}^{\pite,+}$ and terminating when $\norm{\nabla\varphi_{\lambda_k,\mu_k}(\cdot)}{} \leq \tau^\ite$	\label{algo:meta_mininizer}
% 	\State Choose $\lambda^\iite < \lambda^\ite$  \label{algo:lambda_selec}
% 	\State Choose $\mu^\iite < \mu^\ite$ with $\lambda^\iite \geq \mu^\iite$ \label{algo:mu_selec}
% 	\State Choose $\tau^\iite > 0$ such that $\lim_\ite \tau^\ite = 0$ \label{algo:tau_selec}
% 	\State $\ite \leftarrow \iite$
% 	\EndWhile
% \end{algorithmic}
\begin{algorithmic}[1]
\algrenewcommand\algorithmicindent{0.75em}%
\Statex%
	\hspace*{-0.75em}%
	\begin{tabular}[t]{@{}l@{}}%
		Choose $\myvec x^0\in\R^\dimn$ and sequences \(\tau_k,\lambda_k,\mu_k>0\) such that%
	\\
		$\tau_k\searrow0$
		~and~
		$0<\mu_k<\lambda_k\searrow0$
		~with~
		$\lambda^k-\mu^k\searrow0$
	\\
		Set $\ite=1$
	\end{tabular}
\While{ stopping criterion is not satisfied }
	\State\label{algo:meta_mininizer}\vspace*{3pt}%
		\begin{tabular}[t]{@{}l@{}}%
			Starting at $\myvec x^{\pite}$, use a descent method to find
		\\
			$\myvec{x}^{\ite}$ such that $\norm{\nabla\varphi_{\lambda_k,\mu_k}(\myvec x^{\ite})}{}\leq\tau^\ite$
		\end{tabular}
	\State
		$\ite \leftarrow \iite$
\EndWhile
\end{algorithmic}
\end{algorithm}

	\section{Numerical experiments} \label{sec:numer_exps}
		We discuss a number of experiments performed to evaluate the proposed method of \cref{sec:min_algos} in two problems: signal decoding and spectral unmixing.
We make use of the nonconvex functions discussed in \cref{thm:ex_ind__binary,thm:ex_l0} to formulate these problems.
%In our analysis, we focus on metrics of interest for practical purposes, such as reconstruction errors.

\subsection{Signal decoding}
	We approach the problem of binary-signal decoding through a constrained least-squares formulation:
	\begin{equation} \label{eq:signal decoding}
	\underset{\myvec{x} \in \R^{P}}{\text{minimize}} \quad \norm{ \myvec{y} - \myvec{H} \myvec{x}}{}^2,
	\qquad \text{subject to} \quad \myvec{x} \in \{0, 1\}^P, 
	\end{equation}
	where $\myvec{x} \in \R^{P}$,  $\myvec{y} \in \R^{N}$ are the received and transmitted signals, respectively, and $\myvec{H} \in \R^{N \times P}$ is the transfer matrix of a given channel.
	Since this problem's constraint is nonconvex, a convex alternative to this problem is the one given by replacing the constraint with $\myvec{x} \in [0, 1]^P$.
	This alternative generates estimates that are not guaranteed to be binary numbers, so they are usually projected into the original problem's feasible set.
	We randomly generated a number of experiments with varying problem dimensionality (both over- and underdetermined formulations), noise levels, and conditioning of the matrix $\myvec{H}$.
	%, and probability distribution of different bits.
	We compared the proposed method with an \ac{ADMM} implementation, with an \ac{ADMM} implementation of the aforementioned convex relaxation, and with a simple least-squares approach to the problem (i.e., without constraints); we denote the four methods by LL, AN, AR, and LS, respectively.
	We evaluated their performance by computing the \ac{BER} relative to the transmitted signal.
	Overall, we found that the proposed method performed similarly or better than the others when $N < P$, and we present the performance of the three methods over examples of such settings in \cref{tab:BER}.
	In more detail, the problem was set as follows: we started by generating the transmitted signal assuming that it was drawn i.i.d. from a Bernoulli distribution, with equal probability of generating both symbols; we then generated the matrix $\myvec{H}$ assuming that its rows were drawn i.i.d. from a multivariate Gaussian distribution with zero mean and covariance matrix $\myvec{\Sigma} \in \R^{P \times P}$, where $[\myvec{\Sigma}]_{ij} = \rho^{|i-j|}$; finally, we added i.i.d. noise drawn from a Gaussian distribution such that its variance resulted in a given \ac{SNR}.
	We tested the three methods by running the experiment 50 times; at the end, we projected all methods' estimates into the set $\{0, 1\}^P$, i.e., we made $[\myvec{x}]_i^{\text{proj}} = 0$ if $[\myvec{x}]_i^{\text{est}} < 0.5$ and $[\myvec{x}]_i^{\text{proj}} = 1$ if $[\myvec{x}]_i^{\text{est}} \geq 0.5$ for $i = 1, \dots, P$, where $\myvec{x}^{\text{est}}$ denotes a given method's estimate and $\myvec{x}^{\text{proj}}$ denotes its projection.
	The algorithms' parameters were tuned as follows: for LL, we let $\lambda^1 = 10^5$ and $\mu^1 = 0.999 \times \lambda^1$, and reduced these parameters by $90\%$ every iteration; in \cref{sec:min_algos}, for AR and AN, we followed the heuristic proposed in~\cite{Boyd2011} to tune step sizes.
	The initial estimate for AN, AN, and LL was the LS estimate, and we stopped them either if the \ac{RMSE} relative to the transmitted signal was lower than $10^{-9}$ or after $10^3$ iterations.
	The experiments were conducted using MATLAB on a machine with a quad-core Intel Core i7 CPU running at 2.5~GHz and with 16~GB of DDR3L RAM.

	\begin{table}[tbp]
		\scriptsize
		\caption{Average \ac{BER} [\%] for estimates of the four methods.} \label{tab:BER}
		\begin{center}
			\begin{tabular}{ | l | l | l | l ||  c | c | c | c |}
				\hline
				N & P & $\rho$ & \ac{SNR} [dB] & LS & AR & AN & LL \\ \hline
				\multirow{18}{*}{20} & \multirow{6}{*}{40}  & \multirow{3}{*}{0.0} & 30 & 47.15 & \textbf{11.35} & 32.30 & 15.40 \\ 
				 &  &  & 20 & 47.45 & \textbf{14.05} & 33.30 & 17.20 \\ 
				 &  &  & 10 & 47.35 & \textbf{21.90} & 31.30 & \textbf{21.90} \\  \cline{3-8}
				 &  & \multirow{3}{*}{0.5} & 30 & 45.75 & \textbf{12.70} & 32.15 & 16.50 \\ 
				 &  &  & 20 & 46.20 & \textbf{13.95} & 30.85 & 16.75 \\ 
				 &  &  & 10 & 47.90 & 27.80 & 34.20 & \textbf{27.20} \\   \cline{2-8}
				 & \multirow{6}{*}{80}  & \multirow{3}{*}{0.0} & 30 & 50.00 & 37.45 & 40.58 & \textbf{35.95} \\ 
				 &  &  & 20 & 49.02 & 35.95 & 39.42 & \textbf{34.12} \\ 
				 &  &  & 10 & 49.83 & 36.12 & 39.92 & \textbf{35.23} \\   \cline{3-8}
				 &  & \multirow{3}{*}{0.5} & 30 & 48.05 & 34.27 & 38.35 & \textbf{33.25} \\ 
				 &  &  & 20 & 49.48 & 35.73 & 39.67 & \textbf{34.50} \\ 
				 &  &  & 10 & 49.08 & 37.38 & 39.92 & \textbf{37.27} \\   \cline{2-8}
				 & \multirow{6}{*}{100}  & \multirow{3}{*}{0.0} & 30 & 48.84 & 37.86 & 41.26 & \textbf{37.82} \\ 
				 &  &  & 20 & 49.24 & \textbf{39.24} & 42.02 & 39.76 \\ 
				 &  &  & 10 & 49.54 & 40.56 & 43.36 & \textbf{39.98 }\\   \cline{3-8}
				 &  & \multirow{3}{*}{0.5} & 30 & 49.94 & 40.44 & 42.36 & \textbf{40.00} \\ 
				 &  &  & 20 & 49.38 & 39.54 & 41.68 & \textbf{38.28} \\ 
				 &  &  & 10 & 48.60 & 39.20 & 41.56 & \textbf{38.56 }\\
				\hline
			\end{tabular}
			\label{tab:binary_LS}
		\end{center}
	\end{table}

\subsection{Spectral unmixing}
\begin{figure}[tbp]
		\centerline{\includegraphics[width=0.45\textwidth]{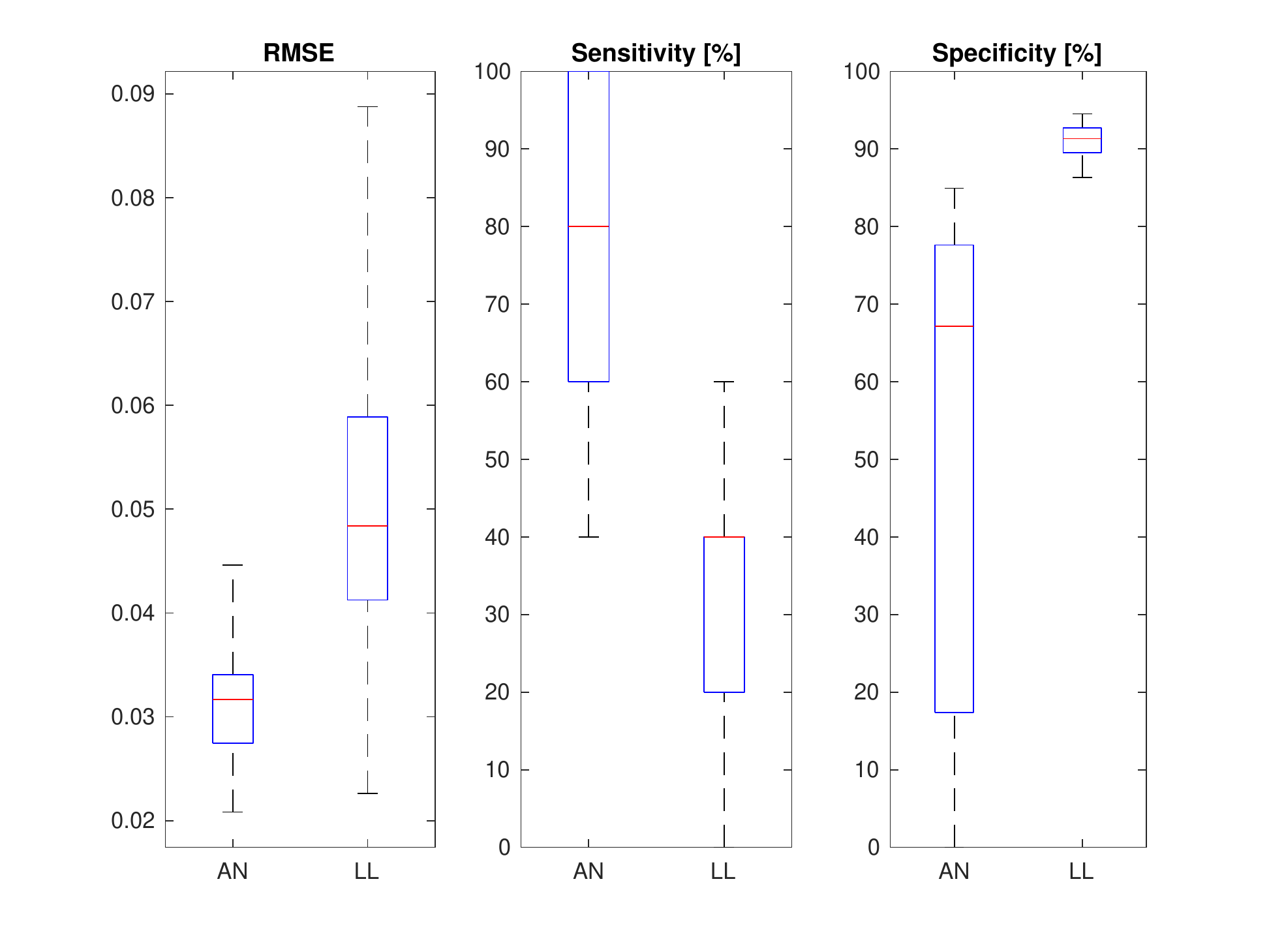}}
		\caption{Average \ac{RMSE}, sensitivity, and specificity.}
		\label{fig:sensitivity_vs_specificity}
	\end{figure}

	Hyperspectral images are multi-channel images with a relatively large number of channels---usually known as spectral bands---corresponding to short frequency ranges along the electromagnetic spectrum.
	Frequently, their spatial resolution is low, and it is of interest to disentangle the different spectral components of a given pixel; a pixel typically corresponds to a mixture of different materials.
	Spectral unmixing techniques produce a set of spectral profiles, one for each material (known as endmember), and a corresponding set of abundances, or percentages of occupation, for each endmember, in each pixel~\cite{Ma2014}.
	We consider that the set of spectral profiles is known through a database of spectral signatures (i.e., a database of reflectance profiles as a function of wavelength), and formulate the spectral unmixing problem pixel-wise as an instance of the Lagrangian formulation of the best--subset-selection problem:
	\begin{equation} \label{eq:l1_unmix}
	\begin{aligned}
	& \underset{\mathbf{a}_j \in \mathbb{R}^{P}}{\text{minimize}} 
	& & \| [\mathbf{Y}_h]_{:j} - \mathbf{U} \mathbf{a}_j \|^2_2 + \beta \| \mathbf{a}_j \|_0,
	\end{aligned}
	\end{equation}
	where $\mathbf{a}_j \in \R^{P}$ is the vector of each endmember's abundances for a given pixel $j$, to be estimated, $\mathbf{U} \in \R^{N \times P}$ is a matrix corresponding to the spectral database, $\mathbf{Y}_h \in \R^{N \times M}$ corresponds to a matrix representation of a hyperspectral image with $M$ pixels (i.e., corresponds to the lexicographical ordering of a 3-D data cube), and $\beta$ is a regularization parameter.
	A typical convex formulation of this problem is as the \ac{LASSO}, where $\| \mathbf{a}_j \|_0$ is replaced by $\| \mathbf{a}_j \|_1$.
	Constraints such as $[\mathbf{a}_j]_i \geq 0$ for $i = 1, \dots, P$ are also adopted but, for the comparison purposes of this work, we ignored them.
	We randomly generated a number of experiments with varying problem dimensionality (both over- and underdetermined formulations) and noise levels.
	We used a real-world spectral dictionary $\mathbf{U}$: a selection of 498 different mineral types from a USGS library, set up as detailed in~\cite{Bioucas2010}.
	We compared the proposed method to an \ac{ADMM} implementation, and we denote the two methods by LL and AN, respectively.
	We evaluated them in two fronts: we computed the \ac{RMSE} relative to the original vector of abundances, and we evaluated the quality of the estimates' support (i.e., the location of its nonzero elements).
	The latter was evaluated by computing the sensitivity and specificity of the estimates as follows: we considered a correctly estimated nonzero as a true positive (TP), and an incorrectly assigned one as a false positive (FP); conversely, a correctly estimated zero was considered as a true negative (TN), and an incorrectly assigned one as a false negative (FN); then, the sensitivity is given by $\frac{TP}{TP+FN}$ and the specificity by $\frac{TN}{FP+TN}$.
	A box plot of the \ac{RMSE}, sensitivity, and specificity for the two methods for one problem setting is given in \cref{fig:sensitivity_vs_specificity}.
	This problem was generated as follows: we start by generating a vector of abundances with $P=224$ and with 5 nonzero elements, where the abundances are drawn from a Dirichlet distribution; we made $N=P$ and added Gaussian noise such that it would result in a \ac{SNR} of 30 dB.
	The parameters were tuned as before except that, for LL, the initial points are $\lambda^1 = 10^3$ and $\mu^1 = \tfrac{\lambda^1}{2}$, and we reduced them by $10\%$ and $82\%$ every iteration, respectively; the initial estimate was considered to be a vector of zeros and the regularization parameter was manually set to $\beta=10^{-6}$, which we found produced estimates with a support similar to the simulated one.
	The proposed method, in virtue of its tuning parameters, allows for some flexibility in how ``aggressively'' it selects nonzero elements.
	In \cref{tab:SS}, and for the problem setting just described, we show how a progressively lower $\lambda^1$---with still making $\mu^1 = \tfrac{\lambda^1}{2}$ and reducing those as indicated---finds estimates with very different sensitivity and specificity.

	\begin{table}[tbp]
		\footnotesize
		\caption{Average \ac{RMSE}, sensitivity [\%], specificity [\%], and the value [$\times 10^{-3}$] of the cost function (CF) for varying $\lambda^1$ in the LL method.}
		\begin{center}
			\begin{tabular}{| l || c | c | c | c | c | c | c | c |}
				\hline
				$\lambda^1$ & $10^4$ & $10^3$ & $10^2$ & $10$ & $1$ & $10^{-1}$ & $10^{-2}$ \\ \hline
				\ac{RMSE} & 0.035 & \textbf{0.032} & 0.034 & 0.047 & 0.056 & 0.044 & 0.038 \\
				Sens. & 38 & 44 & 48 & 54 & 62 & 74 & \textbf{88} \\
				Spec. & \textbf{95.89} &  94.70 & 90.36 &  82.10 &  67.95 & 52.60 & 37.08 \\
				CF & 0.468 & 0.445 & 0.432 & \textbf{0.424} & 0.426 & 0.464 & 0.513 \\
				\hline
				\end{tabular}
				\label{tab:SS}
			\end{center}
	\end{table}

	\section{Conclusions} \label{sec:conclusions}
		We introduced a novel method to address nonsmooth nonconvex composite-minimization problems, based on Lasry--Lions double envelopes.
This method has very broad applications, and it seems suitable for settings of great practical interest in a number of signal processing problems, as suggested by experiments.
When compared to classical methods, the tuning of its inner parameters allows for greater flexibility in the choice of the type of solutions that are sought.
In future work, we will study techniques to automatically tune these parameters.

 	\bibliographystyle{IEEEtran}
	\bibliography{refs}

\acrodef{ADMM}{alternating-direction method of multipliers}
\acrodef{SNR}{signal-to-noise ratio}
\acrodef{i.i.d.}{independent and identically distributed}
\acrodef{RMSE}{root-mean-squared error}
\acrodef{BER}{binary-error rate}
\acrodef{LASSO}{least absolute shrinkage and selection operator}

\end{document}